# Adaptive Certainty-Equivalence Control With Regulation-Triggered Finite-Time Least-Squares Identification, Part II: Analysis

**Iasson Karafyllis[*] and Miroslav Krstic[**]**

[*]Dept. of Mathematics, National Technical University of Athens, Zografou Campus, 15780, Athens, Greece, email: iasonkar@central.ntua.gr

[**]Dept. of Mechanical and Aerospace Eng., University of California, San Diego, La Jolla, CA 92093-0411, U.S.A., email: krstic@ucsd.edu

## Abstract

We present the stability analysis for the new regulation-triggered approach to adaptive control introduced in a companion paper. Due to the fact that the closed-loop system is hybrid, our proofs have essential differences from the conventional adaptive control proofs, where the Lyapunov analysis either encompasses the complete closed-loop state or is done in multiple steps through comparison or Gronwall-Bellman lemmas. In addition, we present a convenient algorithm for checking our parameter-observability assumption, which involves repeated Lie derivatives of appropriate vector fields and can be applied to the class of nonlinear control systems for which at most one unknown parameter appears in each differential equation.

## 1. Introduction

The stability and convergence analysis for the new regulation-triggered adaptive control approach [3] is presented here. The assumptions and the description of the adaptive regulator are briefly restated for convenience in Section 2 and the proofs are given in Section 3. The design innovation in [3] is accompanied with an innovation in the stability analysis here. The proof of the main result in [3], Theorem 3.1, differs essentially from the conventional adaptive control proofs of the past forty years, where Lyapunov analysis is employed, which either encompasses the complete closed-loop state or is done in multiple steps through comparison or Gronwall-Bellman lemmas [5]. We go beyond the distinct stability treatments of the plant state and the parameter estimate that are characteristic of the modular (estimation-based) approaches to adaptive control. The major distinctions come from the fact that, unlike for conventional approaches, the parameter estimate in our approach converges to the true value (under a reasonably unrestrictive parameter-observability condition) even in the absence of persistent excitation, and, just as importantly, does so in finite time, whereas the plant has a conventional, infinite-time regulation property. Another major distinction from most other adaptive approaches is that our controller, and, hence, the closed-loop system, is hybrid, resulting in our devoting a significant part of the proof of our Theorem 3.1 to the closed-loop system's well-posedness and to showing that, apart from $+\infty$, there is no other accumulation point for the times of events. While our general result is on global asymptotic regulation on the plant state, under suitable assumptions we also present results on (local and global) exponential regulation.

In Section 4 (Theorem 4.1) we present a convenient algorithm for checking our parameter-observability assumption, which involves repeated Lie derivatives of appropriate vector fields and can be applied to the class of nonlinear control systems for which at most one unknown parameter appears in each differential equation.

_Notation._

* For a vector $x \in \Re^n$ we denote by $|x|$ its usual Euclidean norm, by $x'$ its transpose.
* $\Re_+$ denotes the set of non-negative real numbers. $Z_+$ denotes the set of non-negative integers.
* We say that a function $V : \Re^n \to \Re_+$ is positive definite if $V(x) > 0$ for all $x \neq 0$ and $V(0) = 0$. We say that a continuous function $V : \Re^n \to \Re_+$ is radially unbounded if the following property holds: "for every $M > 0$ the set $\{x \in \Re^n : V(x) \leq M\}$ is compact". For a given vector field $\Re^n \ni x \to F(x) \in \Re^n$ and a smooth function $V : \Re^n \to \Re$, $L_F V(x)$ denotes the Lie derivative of $V$ along $F$, i.e., $L_F V(x) = \nabla V(x) F(x)$, where $\nabla V(x) = \left( \frac{\partial V}{\partial x_1}(x), \ldots, \frac{\partial V}{\partial x_n}(x) \right)$. The repeated Lie derivative $L_F^{(k)} V(x)$ for certain integer $k > 1$ is defined recursively by the formula $L_F^{(k)} V(x) = L_F L_F^{(k-1)} V(x)$, with $L_F^{(1)} V(x) = L_F V(x)$.
* By $K$ we denote the class of strictly increasing $C^0$ functions $a : \Re_+ \to \Re_+$ with $a(0) = 0$. By $K_\infty$ we denote the class of strictly increasing $C^0$ functions $a : \Re_+ \to \Re_+$ with $a(0) = 0$ and $\lim_{s \to +\infty} a(s) = +\infty$. By $KL$ we denote the set of all continuous functions $\sigma : \Re_+ \times \Re_+ \to \Re_+$ with the properties: (i) for each $t \geq 0$ the mapping $\sigma(\cdot, t)$ is of class $K$ ; (ii) for each $s \geq 0$, the mapping $\sigma(s, \cdot)$ is non-increasing with $\lim_{t \to +\infty} \sigma(s, t) = 0$.

All stability notions used in this paper are the standard stability notions for time-invariant systems (see [4]).

## 2. Assumptions and Description of the Adaptive Regulator

Consider the system

$$\dot{x} = f(x,u) + g(x,u)\theta \qquad x \in \Re^n, u \in \Re^m, \theta \in \Re^l \tag{2.1}$$

where $f : \Re^n \times \Re^m \to \Re^n$, $g : \Re^n \times \Re^m \to \Re^{n \times l}$ are smooth mappings with $f(0,0) = 0$, $g(0,0) = 0$ and $\theta \in \Re^l$ is a vector of constant but unknown parameters.

We first suppose that there exist a smooth mapping $k : \Re^l \times \Re^n \to \Re^m$ with $k(\theta, 0) = 0$ for all $\theta \in \Re^l$ (the nominal feedback controller), two families of continuous, positive definite and radially unbounded functions $V_\theta, Q_\theta : \Re^n \to \Re_+$ parameterized by $\theta \in \Re^l$ (Lyapunov-like functions) with the mappings $\Re^l \times \Re^n \ni (\theta, x) \to Q_\theta(x)$, $\Re^l \times \Re^n \ni (\theta, x) \to V_\theta(x)$ being continuous and such that the following assumptions hold.

**(H1)** *For each $\theta \in \Re^l$, $0 \in \Re^n$ is Globally Asymptotically Stable (GAS) for the closed-loop system*

$$\dot{x} = f(x, k(\theta, x)) + g(x, k(\theta, x))\theta \tag{2.2}$$

*Moreover, for every $\theta \in \Re^l$, $x_0 \in \Re^n$ the solution $x(t) \in \Re^n$ of (2.2) with initial condition $x(0) = x_0$ satisfies the inequality $V_\theta(x(t)) \leq Q_\theta(x_0)$ for all $t \geq 0$.*

**(H2)** *For every non-empty, compact set $\Theta \subset \Re^l$, the following property holds: "for every $M \geq 0$ there exists $R > 0$ such that the implication $V_\theta(x) \leq M, \theta \in \Theta \Rightarrow |x| \leq R$ holds".*

Assumption (H1) is a standard stabilizability assumption (necessary for all possible adaptive control design methodologies). Assumption (H2) is a technical assumption, which requires a "uniform" coercivity property for the family $V_\theta$ on compact sets of $\Re^l$.

In order to be able to estimate the vector of constant but unknown parameters $\theta \in \Re^l$, we need an additional parameter observability assumption.

**(H3)** *There exists a positive integer $N$ such that the following implication holds:*
*"If there exist times $0 = \tau_0 < \tau_1 < ... < \tau_N$, vectors $\theta, d_0, ..., d_N \in \Re^l$ with $d_i \neq 0$ for $i = 0,...,N$ and a right differentiable mapping $x \in C^0\left([0,\tau_N];\Re^n\right) \cap C^1\left([0,\tau_N] \setminus \{\tau_0,...,\tau_N\};\Re^n\right)$ satisfying $\dot{x}(t) = f(x(t), k(\theta + d_i, x(t))) + g(x(t), k(\theta + d_i, x(t)))\theta$ for $t \in [\tau_i, \tau_{i+1})$, $i = 0,...,N-1$, $g(x(t), k(\theta + d_j, x(t)))d_{i+1} = 0$ for all $t \in [\tau_j, \tau_{j+1}]$, $i = 0,...,N-1$, $j = 0,...,i$, then $x(t) = 0$ for all $t \in [0, \tau_N]$."*

Assumption (H3) is an observability assumption for the closed-loop system (2.1) with $u = k(\hat{\theta}, x)$, which guarantees that the only solution for which the vector of constant but unknown parameters $\theta \in \Re^l$ cannot be estimated is the zero solution. Assumption (H3) replaces the well-known "persistency of excitation" condition that is used in many cases for the design of adaptive control schemes.

Next we present the adaptive regulator that is based on assumptions (H1), (H2), (H3). The reasoning behind the construction of the adaptive scheme is described in detail in the companion paper [3].

The control action in the interval between two consecutive events is governed by the nominal feedback $u = k(\theta, x)$ with the unknown $\theta \in \Re^l$ replaced by its estimate $\hat{\theta}$ at the beginning of the interval. Moreover, the estimate $\hat{\theta}$ of the unknown $\theta \in \Re^l$ is kept constant between two consecutive events. In other words, we have

$$\begin{aligned} u(t) &= k\left(\hat{\theta}(\tau_i), x(t)\right) \quad , \quad t \in [\tau_i, \tau_{i+1}), i \in Z_+ \\ \hat{\theta}(t) &= \hat{\theta}(\tau_i) \quad , \quad t \in [\tau_i, \tau_{i+1}), i \in Z_+ \end{aligned} \tag{2.3}$$

where $\{\tau_i \geq 0\}_{i=0}^{\infty}$ is the sequence of the event times satisfying

$$\tau_{i+1} = \min\left(\tau_i + T, r_i\right) \quad , \quad i \in Z_+ \tag{2.4}$$

where $\tau_0 = 0$, $T > 0$ is a positive constant (one of the tunable parameters of the proposed scheme) and $r_i > \tau_i$ is a time instant determined by the event trigger.

Let $a : \Re^n \to \Re_+$ be a continuous, positive definite function (again, one of the tunable parameters of the scheme). The event trigger is described by the equations:

$$r_i := \inf\left\{ t > \tau_i : V_{\hat{\theta}(\tau_i)}(x(t)) = Q_{\hat{\theta}(\tau_i)}(x(\tau_i)) + a(x(\tau_i)) \right\}, \text{ for } x(\tau_i) \neq 0 \tag{2.5}$$

$$r_i := T, \text{ for } x(\tau_i) = 0 \tag{2.6}$$

The description of the event-triggered adaptive control scheme is completed by the parameter update law, which is activated at the times of the events. Let $\tilde{N} > N$ be an (arbitrary; the last of the tunable parameters of the proposed scheme) positive integer that satisfies $\tilde{N} > N$, where $N > 0$ is the positive integer involved in Assumption (H3). Define:

$$p(t,\sigma) := x(t) - x(\sigma) - \int_{\sigma}^{t} f(x(s), u(s))ds \tag{2.7}$$

$$q(t,\sigma) := \int_{\sigma}^{t} g(x(s), u(s))ds \tag{2.8}$$

$$\mu_{i+1} := \min\left\{ \tau_j : j \in \{0,...,i\}, \tau_j \geq \tau_{i+1} - \tilde{N}T \right\} \tag{2.9}$$

$$Z(\tau_{i+1}, \mu_{i+1}) = G(\tau_{i+1}, \mu_{i+1})\,\theta \tag{2.10}$$

where

$$G(\tau_{i+1},\mu_{i+1}) = \int_{\mu_{i+1}}^{\tau_{i+1}} \int_{\mu_{i+1}}^{\tau_{i+1}} q'(t,\sigma)q(t,\sigma)d\sigma\, dt$$

$$Z(\tau_{i+1},\mu_{i+1}) = \int_{\mu_{i+1}}^{\tau_{i+1}} \int_{\mu_{i+1}}^{\tau_{i+1}} q'(t,\sigma)p(t,\sigma)d\sigma\, dt \tag{2.11}$$

The parameter update law is given by the formula

$$\hat{\theta}(\tau_{i+1}) = \arg\min\left\{ \left|\vartheta - \hat{\theta}(\tau_i)\right|^2 : \vartheta \in \Re^l, Z(\tau_{i+1},\mu_{i+1}) = G(\tau_{i+1},\mu_{i+1})\vartheta \right\} \tag{2.12}$$

which, as explained in [3], is a discontinuous, non-recursive, least-squares estimator.

## 3. Analysis of the Adaptive Scheme

In this section, for reader's convenience, we first provide the statements of the results contained in [3].

**Theorem 3.1:** *Consider the control system (2.1) under assumptions (H1), (H2), (H3). Let* $T > 0$ *be a positive constant and let* $a : \Re^n \to \Re_+$ *be a continuous, positive definite function. Finally, let* $\tilde{N} > N$ *be a positive integer that satisfies* $\tilde{N} > N$ *, where* $N > 0$ *is the positive integer involved in Assumption (H3). Then there exists a family of* $KL$ *mappings* $\tilde{\sigma}_{\theta,\hat{\theta}} \in KL$ *parameterized by* $\theta \in \Re^l$*,* $\hat{\theta} \in \Re^l$ *such that for every* $\theta \in \Re^l$*,* $x_0 \in \Re^n$*,* $\hat{\theta}_0 \in \Re^l$ *the solution of the hybrid closed-loop system (2.1) with (2.3), (2.4), (2.5), (2.6), (2.9), (2.12) and initial conditions* $x(0) = x_0$*,* $\hat{\theta}(0) = \hat{\theta}_0$ *is unique, is defined for all* $t \ge 0$ *and satisfies* $|x(t)| \le \tilde{\sigma}_{\theta,\hat{\theta}_0}(|x_0|, t)$ *for all* $t \ge 0$*. Moreover, if* $x_0 \ne 0$ *then* $\hat{\theta}(t) = \theta$ *for all* $t \ge NT$*.*

**Theorem 3.2:** *Consider the control system (2.1) under assumptions (H1), (H2), (H3). Moreover, suppose that for each* $\theta \in \Re^l$*, there exist constants* $M_\theta, \omega_\theta, R_\theta > 0$ *such that for every* $x_0 \in \Re^n$ *with* $|x_0| \le R_\theta$ *the solution of (2.2) with initial condition* $x(0) = x_0$ *satisfies the estimate* $|x(t)| \le M_\theta \exp(-\omega_\theta t)|x_0|$ *for all* $t \ge 0$*; i.e.,* $0 \in \Re^n$ *is Locally Exponentially Stable (LES) for the closed-loop system (2.2). Furthermore, suppose that for every nonempty, compact set* $\Theta \subset \Re^l$ *there exist constants* $R > 0$*,* $K_2 > K_1 > 0$ *such that*

$$K_1|x|^2 \le V_\theta(x) \le Q_\theta(x) \le K_2|x|^2, \quad \text{for all } x \in \Re^n, \theta \in \Theta \text{ with } |x| \le R \tag{3.1}$$

*Let* $T > 0$ *be a positive constant and let* $a : \Re^n \to \Re_+$ *be a continuous, positive definite function that satisfies* $\sup\left\{ |x|^{-2} a(x) : x \in \Re^n, x \ne 0, |x| \le \delta \right\} < +\infty$ *for certain* $\delta > 0$*. Finally, let* $\tilde{N} > N$ *be a positive integer that satisfies* $\tilde{N} > N$ *, where* $N > 0$ *is the positive integer involved in Assumption (H3). Then there exists a family of constants* $\tilde{M}_{\theta,\hat{\theta}}, \tilde{R}_{\theta,\hat{\theta}} > 0$ *parameterized by* $(\theta,\hat{\theta}) \in \Re^l \times \Re^l$*, such that for every* $\theta \in \Re^l$*,* $x_0 \in \Re^n$*,* $\hat{\theta}_0 \in \Re^l$ *with* $|x_0| \le \tilde{R}_{\theta,\hat{\theta}_0}$ *the solution of the hybrid closed-loop system (2.1) with (2.3), (2.4), (2.5), (2.6), (2.9), (2.12) and initial conditions* $x(0) = x_0$*,* $\hat{\theta}(0) = \hat{\theta}_0$ *satisfies the estimate* $|x(t)| \le \tilde{M}_{\theta,\hat{\theta}_0} \exp(-\omega_\theta t)|x_0|$ *for all* $t \ge 0$*.*

**Theorem 3.3:** *Consider system (2.1) under assumptions (H1), (H2), (H3). Moreover, suppose that for each* $\theta \in \Re^l$*,* $0 \in \Re^n$ *is Globally Exponentially Stable (GES) for (2.2) and that for every nonempty, compact set* $\Theta \subset \Re^l$ *there exist constants* $K_2 > K_1 > 0$ *such that*

$$K_1|x|^2 \le V_\theta(x) \le Q_\theta(x) \le K_2|x|^2, \quad \text{for all } x \in \Re^n, \theta \in \Theta \tag{3.2}$$

*Let $T>0$ be a positive constant and let $a:\Re^n \to \Re_+$ be a continuous, positive definite function that satisfies $\sup\left\{ |x|^{-2} a(x): x\in\Re^n, x\neq 0 \right\} < +\infty$. Finally, let $\tilde{N} > N$ be a positive integer that satisfies $\tilde{N} > N$, where $N>0$ is the positive integer involved in Assumption (H3). Then there exists a family of constants $\tilde{M}_{\theta,\hat{\theta}} > 0$ parameterized by $\theta\in\Re^l$, $\hat{\theta}\in\Re^l$, such that for every $\theta\in\Re^l$, $x_0\in\Re^n$, $\hat{\theta}_0\in\Re^l$ the solution of the hybrid closed-loop system (2.1) with (2.3), (2.4), (2.5), (2.6), (2.9), (2.12) and initial conditions $x(0)=x_0$, $\hat{\theta}(0)=\hat{\theta}_0$ satisfies the estimate $|x(t)| \le \tilde{M}_{\theta,\hat{\theta}_0} \exp(-\omega_\theta t)|x_0|$ for all $t\ge 0$.*

**Corollary 3.4:** *Consider the system*

$$\begin{aligned} &\dot{x} = (A+\theta_1 C_1 + ... + \theta_l C_l)x + Bu \\ &x\in\Re^n, u\in\Re^m, \theta=(\theta_1,...,\theta_l)'\in\Re^l \end{aligned} \tag{3.3}$$

*where $A, C_1,....,C_l \in \Re^{n\times n}$, $B\in\Re^{n\times m}$ are constant matrices. Suppose that there exists a family of constants $\omega_\theta > 0$ parameterized by $\theta\in\Re^l$ and a continuous mapping $\Re^l \ni \theta \to M(\theta)\in[1,+\infty)$ such that $\left|\exp\left(t\left(A+\theta_1 C_1+...+\theta_l C_l + BK_\theta\right)\right)\right| \le \exp(-\omega_\theta t)M(\theta)$ for all $t\ge 0$. Moreover, suppose that for every $\theta=(\theta_1,...,\theta_l)'\in\Re^l$, $\hat{\theta}\in\Re^l$, $\vartheta=(\vartheta_1,...,\vartheta_l)'\in\Re^l$ with $\hat{\theta}\neq\theta$ and $\vartheta\neq 0$, the pair of matrices $\left(A+\theta_1 C_1+...+\theta_l C_l + BK_{\hat{\theta}}, \vartheta_1 C_1 + ... + \vartheta_l C_l\right)$ is an observable pair of matrices. Let $a, T>0$ be constants and let $\tilde{N} > 1$ be a positive integer. Let $L:\Re^n \to \Re^{n\times l}$ be the linear operator defined by $L*x = (C_1 x e_1' + ... + C_l x e_l')\in\Re^{n\times l}$ for $x\in\Re^n$ with $e_1' = (1,0,...,0)'\in\Re^l$, ..., $e_l' = (0,...,0,1)'\in\Re^l$. Then there exists a family of constants $\tilde{M}_{\theta,\hat{\theta}} > 0$ parameterized by $\theta\in\Re^l$, $\hat{\theta}\in\Re^l$, such that for every $\theta\in\Re^l$, $x_0\in\Re^n$, $\hat{\theta}_0\in\Re^l$ the solution of the hybrid closed-loop system (3.3) with (2.4), (2.6), (2.9),*

$$\begin{aligned} &u(t) = K_{\hat{\theta}(\tau_i)} x(t) \quad , \quad t\in[\tau_i, \tau_{i+1}), i\in Z_+ \\ &\hat{\theta}(t) = \hat{\theta}(\tau_i) \quad , \quad t\in[\tau_i, \tau_{i+1}), i\in Z_+ \end{aligned} \tag{3.4}$$

$$r_i := \inf\left\{ t>\tau_i : |x(t)| = |x(\tau_i)| \sqrt{a + M^2(\hat{\theta}(\tau_i))} \right\}, \text{ for } x(\tau_i)\neq 0, \tag{3.5}$$

$$\hat{\theta}(\tau_{i+1}) = \arg\min\left\{ \left|\vartheta - \hat{\theta}(\tau_i)\right|^2 : \vartheta\in\Re^l, q(\tau_{i+1},\mu_{i+1}) = Q(\tau_{i+1},\mu_{i+1})\vartheta \right\} \tag{3.6}$$

*where*

$$\begin{aligned} &\dot{z} = x, \ z\in\Re^n \\ &\dot{w} = u, \ w\in\Re^m \\ &y = x - Az - Bw \\ &q(\tau,\mu) = \int_\mu^\tau \int_\mu^\tau \left(L*(z(t)-z(\sigma))\right)' \left(y(t)-y(\sigma)\right) d\sigma\, dt \\ &Q(\tau,\mu) = \int_\mu^\tau \int_\mu^\tau \left(L*(z(t)-z(\sigma))\right)' \left(L*(z(t)-z(\sigma))\right) d\sigma\, dt \end{aligned} \tag{3.7}$$

*and initial conditions $x(0)=x_0$, $\hat{\theta}(0)=\hat{\theta}_0$, $z(0)=0$, $w(0)=0$ satisfies the estimate $|x(t)| \le \tilde{M}_{\theta,\hat{\theta}_0}\exp(-\omega_\theta t)|x_0|$ for all $t\ge 0$. Moreover, if $x_0\neq 0$ then $\hat{\theta}(t)=\theta$ for all $t\ge T$.*

Next, we provide the proofs of all results.

**Proof of Theorem 3.1:** The first claim is a direct consequence of the event trigger given by (2.5) and (2.6). The proof of Claim 1 is straightforward and is omitted.

**Claim 1:** *If a solution $(x(t),\hat{\theta}(t))$ of (2.1) with (2.3), (2.4), (2.5), (2.6), (2.9) and (2.12) is defined on $t\in[0,\tau_i]$ for certain $i\in Z_+$, then the solution is defined on $t\in[0,\tau_{i+1}]$. Moreover, it holds that*

$$V_{\hat{\theta}(\tau_i)}(x(t)) \le Q_{\hat{\theta}(\tau_i)}(x(\tau_i)) + a(x(\tau_i)), \text{ for all } t\in[\tau_i,\tau_{i+1}] \tag{3.8}$$

The next claim clarifies what happens when the parameter estimation error becomes zero at the time of an event.

**Claim 2:** *If a solution* $(x(t),\hat{\theta}(t))$ *of (2.1) with (2.3), (2.4), (2.5), (2.6), (2.9) and (2.12) satisfies* $\hat{\theta}(\tau_i)=\theta$ *for certain* $i\in Z_+$, *then the solution is defined for all* $t\geq 0$ *and satisfies* $\hat{\theta}(t)=\theta$ *for all* $t\geq\tau_i$ *and* $\tau_j=\tau_i+(j-i)T$ *for* $j\geq i$.

**Proof of Claim 2:** Notice that $\hat{\theta}(t)=\theta$ for all $t\in[\tau_i,\tau_{i+1})$. Assume first that $x(\tau_i)\neq 0$. Since for every $\theta\in\Re^l$, $y_0\in\Re^n$ the solution $y(t)\in\Re^n$ of $\dot{y}=f(y,k(\theta,y))+g(y,k(\theta,y))\theta$ with initial condition $y(0)=y_0$ satisfies the inequality $V_\theta(y(t))\leq Q_\theta(y_0)$ for all $t\geq 0$ (recall assumption (H1)), it follows that $V_{\hat{\theta}(\tau_i)}(x(t))\leq Q_{\hat{\theta}(\tau_i)}(x(\tau_i))<Q_{\hat{\theta}(\tau_i)}(x(\tau_i))+a(x(\tau_i))$ for all $t\in[\tau_i,\tau_{i+1}]$. Consequently, it follows from (2.5) that $\tau_{i+1}=\tau_i+T$. The same conclusion follows from (2.6) if $x(\tau_i)=0$. Since equation (2.10) holds, it follows from (2.12) that $\hat{\theta}(\tau_{i+1})=\theta$. Applying induction, we conclude that $\hat{\theta}(t)=\theta$ for all $t\geq\tau_i$ and $\tau_j=\tau_i+(j-i)T$ for all $j\geq i$. $\triangleleft$

The third claim shows that there is actually the time of an event with zero parameter estimation error. This is important because Claim 3 in conjunction with Claim 2 shows that the hybrid system (2.1) with (2.3), (2.4), (2.5), (2.6), (2.9) and (2.12) is a well-posed, forward complete system with no accumulation point for the times of the events other than $+\infty$.

**Claim 3:** *If* $x(0)\neq 0$ *and* $\hat{\theta}(0)\neq\theta$ *then there exists an integer* $i\in\{1,...,N\}$ *that satisfies* $\hat{\theta}(\tau_i)=\theta$.

**Proof of Claim 3:** The proof is made by contradiction. Suppose that $x(0)\neq 0$, $\hat{\theta}(\tau_i)\neq\theta$ for $i=0,1,...,N$. We define

$$d_i:=\hat{\theta}(\tau_i)-\theta,\ \text{for}\ i=0,1,...,N \tag{3.9}$$

It follows from (2.3) and (3.9) that there exist times $0=\tau_0<\tau_1<...<\tau_N$, vectors $\theta,d_0,...,d_N\in\Re^l$ with $d_i\neq 0$ for $i=0,...,N$ and a right differentiable mapping $x\in C^0([0,\tau_N];\Re^n)\cap C^1([0,\tau_N]\backslash\{\tau_0,...,\tau_N\};\Re^n)$ satisfying $\dot{x}(t)=f(x(t),k(\theta+d_i,x(t)))+g(x(t),k(\theta+d_i,x(t)))\theta$ for $t\in[\tau_i,\tau_{i+1})$, $i=0,...,N-1$. We next show that the solution also satisfies $g(x(t),k(\theta+d_j,x(t)))d_{i+1}=0$ for all $t\in[\tau_j,\tau_{j+1}]$, $i=0,...,N-1$, $j=0,...,i$. Assumption (H3) guarantees that $x(t)=0$ for all $t\in[0,\tau_N]$, which contradicts the assumption $x(0)\neq 0$.

Equations (2.4), (2.9) and the fact that $\tilde{N}>N$ implies that $\mu_{i+1}=0$ for $i=0,...,N-1$. Using (2.10), (2.12), we obtain:

$$\left(\int_0^{\tau_{i+1}}\int_0^{\tau_{i+1}}q'(t,\sigma)q(t,\sigma)d\sigma\,dt\right)d_{i+1}=0\ ,\text{for}\ i=0,...,N-1 \tag{3.10}$$

Multiplying (3.10) from the left with $d'_{i+1}$ we get:

$$\int_0^{\tau_{i+1}}\int_0^{\tau_{i+1}}\left|q(t,\sigma)d_{i+1}\right|^2d\sigma\,dt=0\ \text{for}\ i=0,...,N-1 \tag{3.11}$$

Continuity of the mapping $(t,\sigma)\to\left|q(t,\sigma)d_{i+1}\right|^2$ and (3.11) imply that the following equality holds for all $t,\sigma\in[0,\tau_{i+1}]$ and $i=0,...,N-1$:

$$q(t,\sigma)d_{i+1}=0 \tag{3.12}$$

Definition (2.8) in conjunction with (3.12) (which implies that $\frac{d}{dt}q(t,\sigma)d_{i+1}=0$ for all $t\in[0,\tau_{i+1}]\backslash\{0,\tau_1,..,\tau_{i+1}\}$ and $\sigma\in[0,\tau_{i+1}]$, $i=0,...,N-1$) implies that

$$g(x(t),u(t))d_{i+1}=0\ ,\ \text{for all}\ t\in[0,\tau_{i+1}]\backslash\{0,\tau_1,..,\tau_{i+1}\},\ i=0,...,N-1 \tag{3.13}$$

Combining (3.13) with (2.3), (3.9) and exploiting the continuity of the mapping $t \to g(x(t),k(\theta + d_j, x(t)))$, we obtain that the solution satisfies $g(x(t),k(\theta + d_j, x(t)))d_{i+1} = 0$ for all $t \in [\tau_j, \tau_{j+1}]$, $i = 0,...,N-1$, $j = 0,...,i$. $\lhd$

We next notice that (2.10) implies that the parameter update law (2.12) satisfies the following estimate for all $i \in Z_+$:

$$\left|\hat{\theta}(\tau_{i+1}) - \hat{\theta}(\tau_i)\right| \leq \left|\theta - \hat{\theta}(\tau_i)\right| \tag{3.14}$$

Using the triangle inequality and (3.14), we get for $i \in Z_+$:

$$\left|\theta - \hat{\theta}(\tau_{i+1})\right| \leq 2\left|\theta - \hat{\theta}(\tau_i)\right| \tag{3.15}$$

Assumption (H1) and Proposition 7 in [6] implies for every $\theta \in \Re^l$ the existence of functions $a_\theta, \beta_\theta \in K_\infty$ such that for every $x_0 \in \Re^n$ the solution $x(t) \in \Re^n$ of (2.2) with $x(0) = x_0$ satisfies the inequality $|x(t)| \leq a_\theta^{-1}\left(\exp(-t)\beta_\theta\left(|x_0|\right)\right)$ for all $t \geq 0$. It follows from Claim 2 and Claim 3 that for every $\theta, \hat{\theta}(0) \in \Re^l$, $x(0) \neq 0$ there exists an integer $j \in \{0,1,...,N\}$ such that $\hat{\theta}(t) = \theta$ for all $t \geq \tau_j$. Consequently, the solution $(x(t),\hat{\theta}(t))$ of the closed-loop system (system (2.1) with (2.3), (2.4), (2.5), (2.6), (2.9) and (2.12)) with initial condition $x(0) \neq 0$, $\hat{\theta}(0) \in \Re^l$ satisfies

$$|x(t)| \leq a_\theta^{-1}\left(\exp(-(t-\tau_j))\beta_\theta\left(|x(\tau_j)|\right)\right), \text{ for } t \geq \tau_j \tag{3.16}$$

Moreover, applying (3.15) inductively, it follows from Claim 2 and Claim 3 and the fact $j \in \{0,1,...,N\}$, that for every $\theta, \hat{\theta}(0) \in \Re^l$, $x(0) \neq 0$ the solution $(x(t),\hat{\theta}(t))$ of the closed-loop system with initial condition $x(0) \neq 0$, $\hat{\theta}(0) \in \Re^l$ satisfies

$$\left|\hat{\theta}(t) - \theta\right| \leq 2^{N-1}\left|\hat{\theta}(0) - \theta\right|, \text{ for all } t \in [0, \tau_j] \tag{3.17}$$

For $s \geq 0$, $\theta \in \Re^l$, define the following functions for $x \in \Re^n$:

$$\begin{aligned} \tilde{V}_\theta(x;s) &:= \min\left\{V_\vartheta(x) : |\vartheta - \theta| \leq 2^{N-1}s\right\} \\ \tilde{Q}_\theta(x;s) &:= \max\left\{Q_\vartheta(x) : |\vartheta - \theta| \leq 2^{N-1}s\right\} + a(x) \end{aligned} \tag{3.18}$$

Proposition 2.9 on page 21 in [1] implies that the mappings $x \to \tilde{V}_\theta(x;s)$, $x \to \tilde{Q}_\theta(x;s)$ are continuous and positive definite for each fixed $s \geq 0$, $\theta \in \Re^l$. Moreover, assumption (H2) guarantees that for each fixed $s \geq 0$, $\theta \in \Re^l$ the mapping $x \to \tilde{V}_\theta(x;s)$ is radially unbounded. Consequently, Proposition 2.2 on page 107 in [2] implies that for each fixed $s \geq 0$, $\theta \in \Re^l$, there exist functions $\tilde{a}_{\theta,s}, \tilde{\beta}_{\theta,s} \in K_\infty$ such that

$$\tilde{a}_{\theta,s}\left(|x|\right) \leq \tilde{V}_\theta(x;s), \ \tilde{\beta}_{\theta,s}\left(|x|\right) \geq \tilde{Q}_\theta(x;s), \text{ for all } x \in \Re^n \tag{3.19}$$

It follows from (3.8), (3.17), definitions (3.18) and (3.19) that for every $\theta, \hat{\theta}(0) \in \Re^l$, $x(0) \neq 0$ the solution $(x(t),\hat{\theta}(t))$ of the closed-loop system with $x(0) \neq 0$, $\hat{\theta}(0) \in \Re^l$ satisfies the following estimate for all $i = 0,...,j$:

$$\tilde{a}_{\theta,s}\left(|x(t)|\right) \leq \tilde{\beta}_{\theta,s}\left(|x(\tau_i)|\right), \text{ for all } t \in [\tau_i, \tau_{i+1}] \tag{3.20}$$

with $s := \left|\hat{\theta}(0) - \theta\right|$. Define for each fixed $s \geq 0$, $\theta \in \Re^l$ $q_{\theta,s}(r) := \tilde{a}_{\theta,s}^{-1}\left(\tilde{\beta}_{\theta,s}(r)\right)$, $r \geq 0$ (notice that $q_{\theta,s} \in K_\infty$). It follows from (3.20) and the fact $j \in \{0,1,...,N\}$ that the following estimate holds:

$$|x(t)| \leq q_{\theta,s}^{(N)}\left(|x(0)|\right), \text{ for all } t \in [0, \tau_j] \tag{3.21}$$

where $q_{\theta,s}^{(N)} := \underbrace{q_{\theta,s} \circ ... \circ q_{\theta,s}}_{N \text{ times}}$ with $s := \left|\hat{\theta}(0) - \theta\right|$. Combining (3.16), (3.21) the facts $j \in \{0,1,...,N\}$ and $\tau_j \leq \tau_N \leq NT$ (a direct consequence of (2.4)), we obtain that for every $\theta \in \Re^l$, $x(0) \neq 0$, $\hat{\theta}(0) \in \Re^l$ the solution $(x(t),\hat{\theta}(t))$ of the closed-loop system with initial condition $x(0) \neq 0$, $\hat{\theta}(0) \in \Re^l$ satisfies the following estimate for all $t \geq 0$:

$$|x(t)| \le a_\theta^{-1}\left(\exp(-(t-NT))\beta_\theta\left(q_{\theta,s}^{(N)}\left(|x(0)|\right)\right)\right) \tag{3.22}$$

with $s := \left|\hat\theta(0)-\theta\right|$. We notice that estimate (3.22) holds for the case $x(0)=0$ as well. The conclusion of the theorem is a direct consequence of (3.22) and the definition $\sigma_{\theta,\hat\theta}(r,t) := a_\theta^{-1}\left(\exp(-(t-NT))\beta_\theta\left(q_{\theta,s}^{(N)}(r)\right)\right)$ for $r,t\ge 0$ and $s := \left|\hat\theta-\theta\right|$. The proof is complete. ◁

**Proof of Theorem 3.2:** Since all assumptions of Theorem 3.1 hold for Theorem 3.2, all relations in the proof of Theorem 3.1 hold.

Since for each $\theta\in\Re^l$, $0\in\Re^n$ is Locally Exponentially Stable (LES) for (2.2), there exists a family of constants $M_\theta,\omega_\theta,R_\theta>0$ parameterized by $\theta\in\Re^l$, such that for every $\theta\in\Re^l$, $x_0\in\Re^n$ with $|x_0|\le R_\theta$ the solution of (2.2) with $x(0)=x_0$ satisfies the estimate $|x(t)|\le M_\theta\exp(-\omega_\theta t)|x_0|$ for $t\ge 0$. Consequently, the solution $(x(t),\hat\theta(t))$ of (2.1) with (2.3), (2.4), (2.5), (2.6), (2.9) and (2.12) with initial condition $x(0)\ne 0$, $\hat\theta(0)\in\Re^l$ satisfies

$$|x(t)| \le \exp(-\omega_\theta(t-\tau_j))M_\theta\left|x(\tau_j)\right|, \text{ for } t\ge\tau_j \tag{3.23}$$

provided that $\left|x(\tau_j)\right|\le R_\theta$. Define

$$A := \sup\left\{|x|^{-2}a(x) : x\ne 0, |x|\le\delta\right\}. \tag{3.24}$$

Notice that due to (3.1), (3.24) and (3.18) for every $\theta\in\Re^l$ there exist constants $R>0$, $K_2>K_1>0$ such that

$$K_1|x|^2\le\tilde V_\theta(x;s),\ (K_2+A)|x|^2\ge\tilde Q_\theta(x;s),$$
$$\text{for all } x\in\Re^n,\theta\in\Theta \text{ with } |x|\le\min(R,\delta) \tag{3.25}$$

with $s := \left|\hat\theta(0)-\theta\right|$. It follows from (3.8), (3.17), definitions (3.18) and (3.25) that for every $\theta,\hat\theta(0)\in\Re^l$, $x(0)\ne 0$, the solution $(x(t),\hat\theta(t))$ of (2.1) with (2.3), (2.4), (2.5), (2.6), (2.9) and (2.12) with initial condition $x(0)\ne 0$, $\hat\theta(0)\in\Re^l$ satisfies the following estimate for all $i=0,...,j$:

$$|x(t)| \le \sqrt{K_1^{-1}(K_2+A)}\left|x(\tau_i)\right|, \text{ for all } t\in[\tau_i,\tau_{i+1}] \tag{3.26}$$

provided that $\left|x(\tau_i)\right|\le\sqrt{K_1(K_2+A)^{-1}}\min(R,\delta)$. It follows from (3.26) and the fact $j\in\{0,1,...,N\}$ that the following estimate holds for all $t\in[0,\tau_j]$:

$$|x(t)| \le K_1^{-N/2}(K_2+A)^{N/2}|x(0)|\le R_\theta, \tag{3.27}$$

provided that $|x(0)|\le K_1^{N/2}(K_2+A)^{-N/2}\min(R,\delta,R_\theta)$. Combining (3.23), (3.27), the facts $j\in\{0,1,...,N\}$ and $\tau_j\le\tau_N\le NT$ (a direct consequence of (2.4)), we obtain that for every $\theta,\hat\theta(0)\in\Re^l$, $x(0)\ne 0$ the solution $(x(t),\hat\theta(t))$ of (2.1) with (2.3), (2.4), (2.5), (2.6), (2.9) and (2.12) with $x(0)\ne 0$, $\hat\theta(0)\in\Re^l$ satisfies the following estimate:

$$|x(t)| \le \exp(-\omega_\theta t)\exp(\omega_\theta NT)M_\theta K_1^{-N/2}(K_2+A)^{N/2}|x(0)|, \text{ for all } t\ge 0 \tag{3.28}$$

provided that $|x(0)|\le K_1^{N/2}(K_2+A)^{-N/2}\min(R,\delta,R_\theta)$. We notice that estimate (3.28) holds for the case $x(0)=0$ as well. The conclusion of the theorem is a direct consequence of (3.28) and definitions

$$\tilde R_{\theta,\hat\theta} := K_1^{N/2}(K_2+A)^{-N/2}\min(R,\delta,R_\theta),\ \tilde M_{\theta,\hat\theta} := \exp(\omega_\theta NT)M_\theta K_1^{-N/2}(K_2+A)^{N/2}.$$

The proof is complete. ◁

**Proof of Theorem 3.3:** The proof of Theorem 3.3 is almost identical with the proof of Theorem 3.2, except of the fact that no restrictions in the magnitude of $|x|$ are needed for the derivation of all estimates. ◁

**Proof of Corollary 3.4:** System (3.3) is a system of the form (2.1) with $f(x,u)=Ax+Bu$ and $g(x,u)=(L*x)\in\Re^{n\times l}$, where $L:\Re^n\to\Re^{n\times l}$ is the linear operator defined by $L*x=C_1xe_1'+...+C_lxe_l'$ for $x \in \Re^n$ with $e_1' = (1,0,...,0)' \in \Re^l$, ... $e_l' = (0,...,0,1)' \in \Re^l$. We next show that (H1), (H2), (H3) hold for system (3.3).

Indeed, assumptions (H1), (H2) hold with $V_\theta(x) := |x|^2$, $Q_\theta(x) := M^2(\theta)|x|^2$, $k(\theta, x) := K_\theta x$ for $\theta \in \Re^l$, $x \in \Re^n$. This fact is a consequence of the assumptions and the previous definitions. Moreover, assumption (H3) holds with $N = 1$. This follows from our observability assumption, which gives the implication:

"If there exist vectors $\theta \in \Re^l$, $\hat{\theta} \in \Re^l$, $\vartheta = (\vartheta_1,...,\vartheta_l)' \in \Re^l$ with $\hat{\theta} \neq \theta$, $\vartheta \neq 0$, $x_0 \in \Re^n$ and a time $\tau_1 > 0$ for which $(\vartheta_1 C_1 + ... + \vartheta_l C_l)\exp\left(t\left(A + \theta_1 C_1 + ... + \theta_l C_l + BK_{\hat{\theta}}\right)\right)x_0 = 0$ for all $t \in [0, \tau_1]$, then $x_0 = 0$."

The above implication is exactly the implication involved in assumption (H3) with $N=1$. The above definitions guarantee that all assumptions of Theorem 3.3 hold. The conclusions of the corollary are consequences of Theorem 3.1 and Theorem 3.3 applied with $a(x) := a|x|^2$ for $x \in \Re^n$. $\triangleleft$

## 4. An Algorithmic Test of the Parameter-Observability Assumption

As shown in [3], one of the advantages of assumption (H3) is that it can be verified a priori without additional assumptions. Indeed, a convenient algorithmic way of checking assumption (H3) can be given for systems with at most one unknown parameter in each differential equation, i.e., systems of the form

$$\dot{x} = f(x,u) + g_1(x,u)\theta_1 e_{N_1} + ... + g_l(x,u)\theta_l e_{N_l} \tag{4.1}$$

where $e_1 = (1,0,...,0)' \in \Re^n$, ... $e_n = (0,...,0,1)' \in \Re^n$ and the integers $N_i \in \{1,...,l\}$ ($i = 1,...,l$) satisfy the condition

$$i, j \in \{1,...,l\}, \ i \neq j \Rightarrow N_i \neq N_j \tag{4.2}$$

Given arbitrary vectors $\theta, \hat{\theta}_1, \hat{\theta}_2, ..., \hat{\theta}_l \in \Re^l$ and a smooth mapping $k : \Re^l \times \Re^n \to \Re^m$ with $k(\theta, 0) = 0$ for all $\theta \in \Re^l$, we perform the following algorithm which gives us a new vector $z \in \Re^l$. Define $F_z(x) = f(x, k(z,x)) + \sum_{i=1}^{l} g_i(x, k(z,x))\theta_i e_{N_i}$ for $z \in \Re^l$.

<u>Algorithm:</u>

Step 1: *Set* $z = \hat{\theta}_1$. *Let* $I_1 \subseteq \{1,...,l\}$ *be the set of indices* $i \in \{1,...,l\}$ *for which the implication*

$$\left.\begin{array}{l} g_i(x, k(z,x)) = 0 \\ L_{F_z}^{(j)} g_i(x, k(z,x)) = 0 \quad , \quad j = 1,2,... \end{array}\right\} \Rightarrow x = 0 \tag{4.3}$$

*holds.*

Step $s > 1$: *Set* $z = \hat{\theta}_s$. *Then set* $z_i = \theta_i$ *for all* $i \in \bigcup_{p=1}^{s-1} I_p$. *Let* $I_s \subseteq \{1,...,l\}$ *be the set of indices for which (4.3) holds.*

We are in a position to prove the following theorem. Its proof is provided in the following section.

**Theorem 4.1:** *Suppose that for every* $\theta, \hat{\theta}_1, \hat{\theta}_2, ..., \hat{\theta}_l \in \Re^l$ *it holds that* $\bigcup_{p=1}^{l} I_p = \{1,...,l\}$. *Then assumption (H3) holds with* $N = l$ *for system (4.1).*

The following examples show how the above algorithm can be applied to nonlinear systems. The examples also show that the fulfilment of assumption (H3) may or may not impose restrictions on the nominal feedback controller.

**Example 4.2:** Consider the planar linear system

$$\dot{x}_1 = x_2 \quad , \quad \dot{x}_2 = \theta(x_1 + cx_2) + u$$
$$x = (x_1, x_2) \in \Re^2, u \in \Re, \theta \in \Re \tag{4.4}$$

where $c \in \Re$ is a known parameter. System (4.4) is a system with only one unknown parameter ($\theta \in \Re$) appearing in one differential equation and therefore is a system of the form (4.1) for which (4.2) holds. A feedback law that globally exponentially stabilizes the origin for (4.4) is given by

$$k(\theta, x) := -k_1 x_1 - k_2 x_2 - \theta(x_1 + cx_2) \tag{4.5}$$

where $k_1, k_2 > 0$ are constants. Applying the algorithm with $g_1(x,u) = x_1 + cx_2$ and taking repeated Lie derivatives with $u = k(z, x)$, we get in the first step:

$$\left.\begin{aligned} x_1 + cx_2 &= 0 \\ -ck_1 x_1 + (1 - ck_2) x_2 + c(\theta - z)(x_1 + cx_2) &= 0 \end{aligned}\right\}$$
$$\Rightarrow x_1 = -cx_2 \, , (1 - ck_2 + c^2 k_1) x_2 = 0$$

Therefore, implication (4.3) holds provided that $ck_2 - c^2 k_1 \neq 1$. Theorem 4.1 guarantees that (H3) holds for system (4.4) with the feedback law (4.5), provided that $ck_2 - c^2 k_1 \neq 1$. It should be noticed that condition $ck_2 - c^2 k_1 \neq 1$ is a condition on the controller gains and implies that the fulfilment of assumption (H3) may impose restrictions on the nominal feedback controller. ◁

**Example 4.3:** Consider the nonlinear system

$$\dot{x}_1 = x_2 \, , \, \dot{x}_2 = x_1^2 + \theta_1 x_2 + x_3 \, , \, \dot{x}_3 = \theta_2 x_1^2 + u$$
$$x = (x_1, x_2, x_3)' \in \Re^3 \, , \theta = (\theta_1, \theta_2)' \in \Re^2 \, , u \in \Re \tag{4.6}$$

System (4.6) is a system with two unknown parameters ($\theta \in \Re^2$), each one appearing in only one differential equation and therefore is a system of the form (4.1) for which condition (4.2) holds. Applying feedback linearization, we obtain a family of feedback laws, which achieve global asymptotic stabilization and local exponential stabilization of the origin for system (4.6):

$$k(\theta, x) := -k_1 x_1 - k_2 x_2 - 2 x_1 x_2 - (\theta_1 + k_3)(x_1^2 + \theta_1 x_2 + x_3) - \theta_2 x_1^2 \tag{4.7}$$

where $k_1, k_2, k_3 > 0$ are constants with $k_2 k_3 > k_1$. We next apply the algorithm. Indeed, we have $g_1(x,u) = x_2$, $g_2(x,u) = x_1^2$ and we get by taking repeated Lie derivatives with $u = k(z, x)$:

$$\left.\begin{aligned} x_1^2 &= 0 \\ 2x_1 x_2 &= 0 \\ 2x_2^2 + 2x_1(x_1^2 + \theta_1 x_2 + x_3) &= 0 \\ 6x_2(x_1^2 + \theta_1 x_2 + x_3) + 2x_1 L_F(x_1^2 + \theta_1 x_2 + x_3) &= 0 \\ 6(x_1^2 + \theta_1 x_2 + x_3)^2 + 8x_2 L_F(x_1^2 + \theta_1 x_2 + x_3) + 2x_1 L_F^{(2)}(x_1^2 + \theta_1 x_2 + x_3) &= 0 \end{aligned}\right\} \Rightarrow x = 0 \tag{4.8}$$

$$\left.\begin{aligned} x_2 &= 0 \\ x_1^2 + \theta_1 x_2 + x_3 &= 0 \\ \theta_1(x_1^2 + \theta_1 x_2 + x_3) - k_1 x_1 - k_2 x_2 - (z_1 + k_3)(x_1^2 + z_1 x_2 + x_3) + (\theta_2 - z_2) x_1^2 &= 0 \end{aligned}\right\} \Rightarrow k_1 x_1 - (\theta_2 - z_2) x_1^2 = 0, x_2 = 0, x_3 = -x_1^2 \tag{4.9}$$

It follows from (4.8) that $2 \in I_1$. However, (4.9) does not necessarily imply that $x = 0$ and we conclude that $I_1 = \{2\}$ if $z_2 \neq \theta_2$. Continuing with the second step of the algorithm and taking $z_2 = \theta_2$, we notice that (4.9) implies $x = 0$ and consequently $I_1 \cup I_2 = \{1\}$. Hence, Theorem 4.1 guarantees that assumption (H3) holds for (4.6) with the feedback law given by (4.7). Notice that for this system the fulfilment of assumption (H3) does not impose any restriction on the nominal feedback controller given by (4.7). ◁

We next provide the proof of Theorem 4.1.

**Proof of Theorem 4.1:** Due to (4.1), (4.2), the verification of assumption (H3) with $N = l$ requires to show the following implication:

"If there exist times $0 = \tau_0 < \tau_1 < ... < \tau_l$, vectors $\theta = (\theta_1,...,\theta_l)'$, $d_0 = (d_{0,1},...,d_{0,l})'$, …, $d_l = (d_{l,1},...,d_{l,l})' \in \Re^l$ with $d_i \neq 0$ for $i = 0,...,l$ and a right differentiable mapping $x \in C^0\left([0,\tau_l];\Re^n\right) \cap C^1\left([0,\tau_l] \setminus \{\tau_0,...,\tau_l\};\Re^n\right)$ satisfying $\dot{x}(t) = f(x(t),k(\theta + d_i, x(t))) + \sum_{p=1}^{l} g_p(x(t),k(\theta + d_i, x(t)))\theta_p e_{N_p}$ for $t \in [\tau_i,\tau_{i+1})$, $i = 0,...,l-1$, $g_p(x(t),k(\theta + d_j, x(t)))d_{i+1,p} = 0$ for all $t \in [\tau_j,\tau_{j+1}]$, $p = 1,...,l$, $i = 0,...,l-1$, $j = 0,...,i$, then $x(t) = 0$ for all $t \in [0,\tau_l]$."

Suppose that there exist times $0 = \tau_0 < \tau_1 < ... < \tau_l$, vectors $\theta = (\theta_1,...,\theta_l)'$, $d_0 = (d_{0,1},...,d_{0,l})'$, …, $d_l = (d_{l,1},...,d_{l,l})' \in \Re^l$ with $d_i \neq 0$ for $i = 0,...,l$ and a right differentiable mapping $x \in C^0\left([0,\tau_l];\Re^n\right) \cap C^1\left([0,\tau_l] \setminus \{\tau_0,...,\tau_l\};\Re^n\right)$ satisfying $\dot{x}(t) = f(x(t),k(\theta + d_i, x(t))) + \sum_{p=1}^{l} g_p(x(t),k(\theta + d_i, x(t)))\theta_p e_{N_p}$ for $t \in [\tau_i,\tau_{i+1})$, $i = 0,...,l-1$, $g_p(x(t),k(\theta + d_j, x(t)))d_{i+1,p} = 0$ for all $t \in [\tau_j,\tau_{j+1}]$, $p = 1,...,l$, $i = 0,...,l-1$, $j = 0,...,i$. Moreover, suppose that the application of the algorithm with $\theta \in \Re^l$ and $\hat{\theta}_1 = \theta + d_0$,…, $\hat{\theta}_l = \theta + d_{l-1}$ gives $\bigcup_{p=1}^{l} I_p = \{1,...,l\}$.

Step 1 of the algorithm with $\hat{\theta}_1 = \theta + d_0$, in conjunction with the fact that $g_p(x(t),k(\theta + d_0, x(t)))d_{i,p} = 0$ for all $t \in [0,\tau_1]$, $p \in I_1$, $i = 1,...,l$ and implication (4.3) guarantees that $x(t) = 0$ for all $t \in [0,\tau_l]$ when $d_{i,p} \neq 0$. Therefore, next we consider the case $I_1 \neq \{1,...,l\}$, $d_{i,p} = 0$ for $p \in I_1$, $i = 1,...,l$.

Step 2 of the algorithm with $\hat{\theta}_2 = \theta + d_1$, in conjunction with the facts that $d_{i,p} = 0$ for $p \in I_1$, $i = 1,...,l$, $g_p(x(t),k(\theta + d_1, x(t)))d_{i,p} = 0$ for all $t \in [\tau_1,\tau_2]$, $p \in I_2$, $i = 2,...,l$ and implication (4.3) guarantees that $x(t) = 0$ for all $t \in [0,\tau_l]$ when $d_{i,p} \neq 0$. Therefore, next we consider the case $I_1 \cup I_2 \neq \{1,...,l\}$, $d_{i,p} = 0$ for all $p \in I_1 \cup I_2$, $i = 2,...,l$.

Continuing in this way, at the final step we conclude that $x(t) = 0$ for all $t \in [0,\tau_l]$, because we cannot have $d_{l,p} = 0$ for all $p \in \bigcup_{q=1}^{l} I_q = \{1,...,l\}$. The proof is complete. $\triangleleft$

## 5. Concluding Remarks

We present the proofs of state regulation and parameter convergence for our new regulation-triggered approach to adaptive nonlinear control [3], as well as a Lie derivative-based algorithm for checking our parameter-observability assumption. Future work may address the relaxation of assumption (H3) as well as the development algorithms for checking the parameter observability assumption for nonlinear systems with more than one parameter per differential equation.

**Acknowledgments:** The authors would like to thank Professor John Tsinias for the useful exchange of ideas and for suggesting Corollary 3.4.